\documentclass[a4paper, 11pt]{amsart}                

\usepackage{caption}
\usepackage{graphicx}
\usepackage{amsmath}
\usepackage{amsfonts}
\usepackage{amssymb}
\usepackage{color}
\usepackage{mathrsfs}
\usepackage{graphicx}
\usepackage{epsfig}
\usepackage[active]{srcltx}
\usepackage{enumerate}
\usepackage{url}

\newtheorem{theorem}{Theorem}

\newtheorem{lemma}[theorem]{Lemma}

\newtheorem{example}[theorem]{Example}
\newtheorem{corollary}[theorem]{Corollary}

\begin{document}

\title{On a class of robust nonconvex quadratic optimization problems
}

\author{F. Flores-Baz\'an$^{1}$ \and Y. Garc\'{\i}a$^{2}$ \and A. P\'erez$^{3}$}

\address{$^{1,3}$Departamento de Ingenier\'{i}a Matem\'atica, Universidad de Concepci\'on, Casilla 160-C,
Concepci\'on, Chile }

\address{$^{2}$ Universidad del Pac\'{\i}fico, Jir\'on S\'anchez Cerro 2050, Jes\'us Mar\'{\i}a, Lima, Per\'u (2)}
\email{fflores@ing-mat.udec.cl}  
 \email{arielperez@udec.cl}          
  \email{garcia\_yv@up.edu.pe}

\maketitle

\begin{abstract}
Let us consider the following robust nonconvex quadratic optimization problem:
\begin{equation*}
\begin{split}
\min &~ \dfrac{1}{2} x^\top Ax+a^\top x  \\
\text{s.t.}~ & \alpha\leq\dfrac{1}{2}x^\top (B_1+\mu B_2)x+(b_1+\delta b_2)^\top x  \leq\beta,~ \forall~ \mu\in [\mu_1,\mu_2],\forall~\delta\in[\delta_1,\delta_2],
\end{split}
\end{equation*}
where $A$, $B_1$, $B_2$ are real symmetric matrices, $\mu_1,\mu_2,\delta_1,\delta_2,\alpha$, $\beta\in\mathbb{R}$ satisfying
$\mu_1\leq \mu_2$, $\delta_1\leq\delta_2$ and $\alpha<\beta$. We 
establish the robust alternative result; the robust S-lemma and the robust optimality for the above nonconvex problem.
\end{abstract}

\def \Pr{\operatornamewithlimits{Pr}}
\def \Min{\operatornamewithlimits{Min}}
\def \Pr{\operatornamewithlimits{Pr}}
\def \l{1\! \!{\rm l}}
\def \minm{\operatornamewithlimits{min}}
\def \lMin{\operatornamewithlimits{\emph{l-}Min}}
\def \lMax{\operatornamewithlimits{\emph{l-}Max}}
\def \DMin{\operatornamewithlimits{DMin}}
\def \lDMin{\operatornamewithlimits{\emph{l-}DMin}}
\def \lDMax{\operatornamewithlimits{\emph{l-}DMax}}
\def \Max{\operatornamewithlimits{Max}}
\def \maxm{\operatornamewithlimits{max}}
\def \DMax{\operatornamewithlimits{DMax}}
\def \Minl{\operatornamewithlimits{Minl}}
\def \DMinl{\operatornamewithlimits{DMinl}}
\def \uMin{\operatornamewithlimits{\emph{u-}Min}}
\def \co{\operatornamewithlimits{conv}}
\def \cono{\operatornamewithlimits{cono}}
\def \Dom{\operatornamewithlimits{Dom}}
\def \Graf{\operatornamewithlimits{Graf}}
\def \dom{\operatornamewithlimits{dom}}
\def \Epi{\operatornamewithlimits{Epi}}
\def \convex{\operatornamewithlimits{convex}}
\def \inte{\operatornamewithlimits{int}}
\def \cl{\operatornamewithlimits{cl}}
\def \diam{\operatornamewithlimits{diam}}
\def \infi{\operatornamewithlimits{2_inf}}
\def \Inf{\operatornamewithlimits{Inf}}
\def \Limsup{\operatorname{Limsup}}
\def \Liminf{\operatorname{Liminf}}
\def \limi{\operatorname{lim}}
\def \Rec{\operatornamewithlimits{Rec}}
\def \Ker{\operatornamewithlimits{Ker}}
\def \cone{\operatornamewithlimits{cone}}

\section{Introduction and basic notation} \label{secc01}
Robust optimization arises as a deterministic approach when addressing an optimization problem under uncertainty data. This paper revisites 
the following robust optimization problem:
\begin{equation}\label{prob:00}
\min\Big\{\dfrac{1}{2}x^\top Ax+a^\top x :~\alpha\leq
\dfrac{1}{2}x^\top Bx+b^\top x\leq\beta,~\forall~(B,b)\in{\mathcal B}_b\Big\},
\end{equation}
where $\mathcal{B}_b\doteq\{B_1+\mu B_2:\mu\in[\mu_1,\mu_2]\}\times\{b_1+\delta b_2:\delta\in[\delta_1,\delta_2]\}$, with all the matrices being real symmetric, $a,b\in\mathbb{R}^n$ and 
$\alpha,\beta,\delta_1,\delta_2,\mu_1,\mu_2$ are given real numbers. 
Optimization problems that can be  modeled by quadratic functions appear, for instance, in \cite{Beck,mat,slemma,Polyak,sid}.

Problem \eqref{prob:00} includes that examined in \cite{jeya2013}:
\begin{equation}\label{prob:01}
\min\Big\{\dfrac{1}{2}x^\top Ax+a^\top x :~
\dfrac{1}{2}x^\top Bx+b^\top x\leq\beta,~\forall~(B,b)\in{\mathcal B}_b\Big\}.
\end{equation}
Theorem 5.1 in \cite{jeya2013} provides a characterization of robust optimality for
problem \eqref{prob:01} under the convexity of the set
 $$
\Big\{(x^\top H_0x,x^\top H_{1}x,x^\top H_{2}x):x\in\mathbb{R}^{n+1}\Big\},
$$
where 
$$
H_0=\begin{pmatrix}A & a\\
a^\top  &  2\gamma \end{pmatrix},
~
H_{1}=\begin{pmatrix}B_1+\mu_1 B_2 & b_1+\delta_1 b_2\\
(b_1+\delta_1b_2)^\top  &  -2\beta \end{pmatrix},~
H_{2}=\begin{pmatrix}B_1+\mu_2 B_2 & b_1+\delta_2 b_2\\
(b_1+\delta_2b_2)^\top  &  -2\beta \end{pmatrix}
$$
with $\gamma=-f(\overline{x})$. Here, $\overline{x}$ is a feasible point of the robust optimization problem \eqref{prob:01}, which is either to be supposed becoming an optimal solution, or to be optimal for deriving optimality conditions.

Certainly the presence of the matrices $H_0, H_1, H_2$ is because the authors in \cite{jeya2013}  homogenize problem \eqref{prob:01} in order to apply  the Dines convexity theorem (\cite{Dines})  valid for quadratic forms. Finally, we realize there is a gap in the proof of Theorem 5.1 in \cite{jeya2013}, but we were unable to find a counterexample to such a result under their assumptions. This is discussed in detail after Theorem \ref{minimocorregido} in Section \ref{secc:uni}. Notice that the approach employed in  \cite{jeya2013} was also applied in \cite{nuevo}.

We have to point out that problem \eqref{prob:00} (and so problem \eqref{prob:01}) was studied without homogenizing the problem thanks to the convexity result established in \cite[Theorem 4.19]{Flores} valid for inhomogeneous quadratic functions. This  allows us to impose the convexity of
a set being the image of $\mathbb{R}^n$ via the  inhomogeneous quadratic functions.

Associated to problem \eqref{prob:00}, purposes of the present paper are to establish an alternative robust result (Theorem \ref{teoralt}), a  robust S-lemma (Theorem \ref{slema}) and a characterization of robust optimality (Theorem \ref{minimonohomoalt}), for problem \eqref{prob:00}. Finally, we provide a counterexample (Example \ref{ex100}) to the argument employed in the proof of Theorem 5.1 in \cite{jeya2013} related to problem \eqref{prob:01}.

Thus, the structure of the present paper is as follows. Section \ref{sect02} establishes the convexity of images for quadratic mappings by applying the Ramana-Goldman criterion \cite{rg1995} (see also \cite[Theorem 2.1 ]{Flores2}). The main results are presented in  Section \ref{sect003}, and Section \ref{secc:uni} revisites problem 
\eqref{prob:01} discussed in \cite{jeya2013}.

\section{A preliminary result: convexity of images}\label{sect02}

By ${\mathcal S}^n$ we denote the set of symmetric matrices of order $n\in\mathbb{N}$ with real entries; ${\mathcal S}_+^n$ denotes the subset of 
${\mathcal S}^n$ whose elements are positive semidefinite matrices, and we write
$A\succeq 0$ if $A\in {\mathcal S}_+^n$; and ${\mathcal S}_{++}^n$ stands for the matrices in ${\mathcal S}^n$ that are positive definite, and in this case we write
$A\succ 0$ if $A\in {\mathcal S}_{++}^n$.

It is our purpose to prove the convexity of images for quadratic mappings under the Ramana-Goldman criterion \cite{rg1995} (see also \cite[Theorem 2.1 ]{Flores2}). To that end, we are 
given $A_i\in{\mathcal S}^n$, $b_i\in\mathbb{R}^n$, $c_i\in\mathbb{R}$ for $i=0, 1,\ldots,m$, we set 
$$
M_i=\begin{pmatrix}A_i & b_i\\
b_i^\top  &  2c_i\end{pmatrix}, f_i(x)= x^\top A_i x + 2b_i^\top x,\; \overline{f}_i(x)= x^\top A_i x,~x\in\mathbb{R}^n.
$$
Furthermore, let us consider the function
$$
G(x,t)= (g_0(x,t),  g_1(x,t),  \ldots,  g_m(x,t)), ~~(x,t)\in\mathbb{R}^n\times\mathbb{R}, 
$$
where $g_i$ is defined by 
$g_i(x,t)= \begin{pmatrix} x\\t \end{pmatrix}^\top M_i\begin{pmatrix} x\\t \end{pmatrix} $.

\begin{lemma} \label{convex_R_5} Let $A_i\in {\mathcal S}^n$,  $c_i \in \mathbb{R}$, 
$b_i\in\mathbb{R}^n$,  $i= 0, 1, \dots, m$ be as above. Set 
$$
F(x)= (f_0(x),  f_1(x),  \ldots,  f_m(x)),~~\overline{F}(x)= (\overline{f}_0(x),  \overline{f}_{1}(x), \ldots,  \overline{f}_{m}(x)).
$$
If  $F(\mathbb{R}^n)$ and  $\overline{F}(\mathbb{R}^n)$ are convex then  $G(\mathbb{R}^{n+1})$ is convex.
\end{lemma}
\begin{proof} Set
	 $\Lambda:= F(\mathbb{R}^n)$ and  $\overline{\Lambda}:= \overline{F}(\mathbb{R}^n)$ and $\Omega:= G(\mathbb{R}^{n+1})$.	
	
	As $\Lambda$ is convex, by the convexity criterion due to Ramana-Goldman  (see also \cite[Theorem 2.1 ]{Flores2}), 
	$\Lambda+\overline{\Lambda}= \Lambda$.  We easily get that for $(x,t) \in \mathbb{R}^{n+1}$,
	\[
	\begin{pmatrix} x\\t \end{pmatrix}^\top M_i\begin{pmatrix} x\\t \end{pmatrix} = x^\top A_ix+2tb_i^{\top} x+2c_i t^2,~~i= 0,1, \dots, m.
	\]
	By setting $\overline{\gamma} = 2(c_0, c_1, \dots,  c_m)$, we obtain
	\begin{equation}\label{contain}
	\Lambda +\overline{\gamma} \subseteq \Omega.
	\end{equation}
	Let ${z}_1 = G(x_1, t_1)$,  ${z}_2 = G(x_2, t_2)$ be any elements in $\Omega$ and let $\lambda \in ~]0,1[$. We distinguish three cases.  \\
	$(i)$: $t_1 \neq 0$ and $t_2 \neq 0$. Then,
	$
	{z}_1= t_1^2(F(x_1/t_1)+ \overline{\gamma})
	$
	and 
	$
	{z}_2= t_2^2(F(x_2/t_2)+ \overline{\gamma}).
	$\\
	The convexity of  $\Lambda$ implies that  
	
	\begin{align*}
	\frac{\lambda}{\lambda t_1^2 +(1-\lambda)t_2^2} &{z}_1+ \frac{1-\lambda}{\lambda t_1^2 +(1-\lambda) t_2^2} {z}_2 \\
	& =  
	\frac{\lambda t_1^2}{\lambda t_1^2 +(1-\lambda)t_2^2}F(x_1/t_1)+\frac{(1-\lambda)t_2^2}{\lambda t_1^2 +(1-\lambda)t_2^2}F(x_2/t_2) + \overline{\gamma} \in  \Lambda + \overline{\gamma}.
	\end{align*}
	Taking into account \eqref{contain} and the fact the $\Omega$ is a cone,  we obtain
	\[
	\lambda {z}_1+ (1-\lambda){z}_2 \in \mathbb{R}_{++}( \Lambda +\overline{\gamma})  \subseteq \Omega.
	\]
	$(ii)$: $t_1=t_2= 0$.  Then,  ${z}_1, {z}_2 \in \overline{\Lambda}$, and because of the convexity of $\overline{\Lambda}$,  we get
	\[
	\lambda \left({z}_1+(1-\lambda){z}_2\right)  \in \overline{\Lambda} \subseteq \Omega.
	\]
	$(iii)$: $t_1\neq 0$ and $t_2=0$. Then,  since $\overline{\Lambda}$ is a cone,
	\[
	\lambda {z}_1+(1-\lambda) {z}_2 \in   \lambda t_1^2( \Lambda +\overline{\gamma}) +  \lambda t_1^2 \overline{ \Lambda}  \subseteq \lambda t_1^2( \Lambda +\overline{\gamma}) \subseteq \Omega.
	\]
	This completes the proof  that $\Omega$ is convex.
	\qed
\end{proof}

Part $(a)$ of the following result is exactly Theorem 2.3 $(i)$ in \cite{Beck}, and $(b)$ is a consequence of the previous lemma.

\begin{corollary} \label{cor_R_5}
	Let the same  hypotheses  of Lemma \ref{convex_R_5} be satisfied.  Let  $\rho_i\in \mathbb{R}$,  for $i=1, \dots, m$.
	If $n\geq m+1$,  $A_0\in {\mathcal S}_{++}^n$,  $A_i= \rho_i A_0$ for $i=1, \dots, m$,  then
	\begin{itemize}
		\item[$(a)$] $F(\mathbb{R}^n)$ and  $\overline{F}(\mathbb{R}^n)$ are convex.
		\item[$(b)$] $G(\mathbb{R}^{n+1})$ is convex.
	\end{itemize}
\end{corollary}
\begin{proof} By assumption on $A_i$, we can apply \cite[Theorem 2.3 (i)]{Beck} to obtain the convexity of $F(\mathbb{R}^n)$ and  $\overline{F}(\mathbb{R}^n)$. Then, $(b)$ follows from Lemma \ref{convex_R_5}.
\qed
\end{proof}

\section{The main results}\label{sect003}

Denote the function:
$$
f(x)\doteq\dfrac{1}{2}x^\top Ax+a^\top x
$$
and let us define the following  matrices in ${\mathcal S}^{n+1}$:

$$
H_0\doteq\begin{pmatrix}A & a\\
a^\top  &  2\gamma \end{pmatrix},~
W(\delta,\lambda)\doteq\begin{pmatrix}B_1 & b_1+\delta b_2\\
(b_1+\delta b_2)^\top  &  -2\lambda \end{pmatrix},~
W_2\doteq\begin{pmatrix}B_2 & 0\\
0^\top  &  0 \end{pmatrix},
$$
and set
$$
W_{1_{\beta}}=W(\delta_1,\beta);~
W_{2_{\beta}}=W(\delta_2,\beta);~ 
W_{1_{\alpha}}=W(\delta_1,\alpha);~
W_{2_{\alpha}}=W(\delta_2,\alpha).
$$

The following set will play an important role in the following.
$$
\Omega_W\doteq\Big\{\Big(\dfrac{1}{2}y^\top H_0y
,\max_{\mu\in[\mu_1,\mu_2]}\dfrac{1}{2}y^\top (W_{1_{\beta}}+\mu W_2)y),\max_{\mu\in[\mu_1,\mu_2]}\dfrac{1}{2}y^\top (W_{2_{\beta}}+\mu W_2)y,
$$
\begin{equation}
-\min_{\mu\in[\mu_1,\mu_2]}\dfrac{1}{2}y^\top (W_{1_{\alpha}}+\mu W_2)y,-\min_{\mu\in[\mu_1,\mu_2]}\dfrac{1}{2}y^\top (W_{2_{\alpha}}+\mu W_2)y\Big):
y\in\mathbb{R}^{n+1}\Big\}+\rm{int} \ \mathbb{R}^5_+.
\end{equation}
By Corollary 1 in \cite{fbp2022}, $\Omega_W$ is convex if the set
$$\Omega_\mu\doteq
\Big\{\Big(y^\top H_0y,y^\top(W_{1_{\beta}}+\mu_1 W_2)y,
y^\top(W_{2_{\beta}}+\mu_1 W_2)y, y^\top (W_{1_{\beta}}+\mu_2 W_2)y,
$$
$$
y^\top (W_{2_{\beta}}+\mu_2 W_2)y,
-y^\top(W_{1_{\alpha}}+\mu_1 W_2)y,
-y^\top(W_{2_{\alpha}}+\mu_1 W_2)y, -y^\top (W_{1_{\alpha}}+\mu_2 W_2)y,
$$
$$
-y^\top (W_{2_{\alpha}}+\mu_2 W_2)y\Big):y\in\mathbb{R}^{n+1} \Big\}+\rm{int} \ \mathbb{R}^9_+
$$
is so.

\begin{theorem} $($A robust alternative result$)$\label{teoralt} 
 Let $A, B_1,B_2\in {\mathcal S}^n$, $a, b_1, b_2\in\mathbb{R}^n$ and $\gamma, \alpha, \beta, \mu_1, \mu_2,\delta_1,\delta_2\in\mathbb{R}$, with $\mu_1\leq\mu_2$, $\delta_1\leq\delta_2$ and $\alpha<\beta$. Assume that $\Omega_W$ is convex. Then, exactly one of the two following assertions hold:
	\begin{itemize}
		\item[$(a)$] $\exists~ x\in\mathbb{R}^n: \dfrac{1}{2}x^\top Ax+a^\top x + \gamma<0$, $\alpha<\dfrac{1}{2}x^\top (B_1+\mu B_2)x+(b_1+\delta b_2)^\top x<\beta$, $\forall~\mu\in [\mu_1,\mu_2]$, $\forall~\delta\in [\delta_1,\delta_2]$.
		\item[$(b)$] $\exists~(\lambda_0,\lambda_1,\lambda_2)\in\mathbb{R}^3_+\backslash\{0\}, \exists~\mu_{\alpha},\mu_{\beta}\in [\mu_1,\mu_2 ], \exists~\delta_{\alpha},\delta_{\beta}\in [\delta_1,\delta_2 ]:~ \forall~ x\in \mathbb{R}^n$
		$$
		\lambda_0\Big(\dfrac{1}{2}x^\top Ax+a^\top x+\gamma\Big) +\lambda_1\Big(	\dfrac{1}{2}x^\top (B_1+\mu_{\beta} B_2)x+(b_1+\delta_{\beta} b_2)^\top x-\beta\Big)+
		$$
		$$ \lambda_2\Big(\alpha-\Big(\dfrac{1}{2}x^\top (B_1+\mu_{\alpha} B_2)x+(b_1+\delta_{\alpha} b_2)^\top x \Big)\Big)\geq 0,
		$$
	\end{itemize}
	where $\mu_{\alpha}+\mu_{\beta}=\mu_1+\mu_2$.\\
	\\
	In addition, we observe that $(b)$ may be written equivalently
	as
	$$
	\lambda_0 A+\lambda_1(B_1+\mu_{\beta}B_2)-\lambda_2(B_1+\mu_{\alpha}B_2)\succeq 0 \ \rm{and}$$
	$$ \exists~ \overline{x}\in\mathbb{R}^n: \Big( \lambda_0A+\lambda_1(B_1+\mu_{\beta}B_2)-\lambda_2(B_1+\mu_{\alpha}B_2)\Big)\overline{x}+\lambda_0 a + \lambda_1(b_1+\delta_{\beta}b_2)+\lambda_2(b_1+\delta_{\alpha} b_2)=0.$$
\end{theorem}
\begin{proof} It is obvious that both statements $(a)$ and 
	$(b)$ cannot be fulfilled simultaneously. Thus, we must
	check that if $(a)$ does not hold, $(b)$ does.\\
	\\
	\textbf{Step 1: The homogenization system.} If $(a)$ does not hold, then there exists no  $x\in\mathbb{R}^n$ such that for all $\mu\in[\mu_1,\mu_2]$ and all $\delta\in[\delta_1,\delta_2]$
	$$
	\dfrac{1}{2}x^\top Ax+a^\top x + \gamma<0, \quad \dfrac{1}{2}x^\top (B_1+\mu B_2)x+(b_1+\delta b_2)^\top x<\beta,
	$$
	$$
	-\left(\dfrac{1}{2}x^\top (B_1+\mu B_2)x+(b_1+\delta b_2)^\top x\right)<-\alpha.
	$$
	By setting $\mathcal{B}_b\doteq\{B_1+\mu B_2:\mu\in[\mu_1,\mu_2]\}\times\{b_1+\delta b_2:\delta\in[\delta_1,\delta_2]\}$, the previous is equivalent to the nonexistence of $x\in\mathbb{R}^n$ such that
	$$
	\dfrac{1}{2}x^\top Ax+a^\top x + \gamma<0,~~	\max\left\{\dfrac{1}{2}x^\top Bx+b^\top x-\beta: (B,b)\in\mathcal{B}_b\right\}<0,
	$$
	$$
	-\min\left\{\dfrac{1}{2}x^\top Bx+b^\top x-\alpha: (B,b)\in\mathcal{B}_b\right\}<0.
	$$
	We claim that the following homogeneous system in $\mathbb{R}^{n+1}$:	
	\begin{equation}\label{sis:hom1}
	\dfrac{1}{2}x^\top Ax+ta^\top x + t^2\gamma<0,~ \max\left\{\dfrac{1}{2}x^\top Bx+tb^\top x-t^2\beta:(B,b)\in\mathcal{B}_b\right\}<0,
	\end{equation}
	\begin{equation}\label{sis:hom2}
	-\min\left\{\dfrac{1}{2}x^\top Bx+tb^\top x-t^2\alpha: (B,b)\in\mathcal{B}_b\right\}<0
	\end{equation}
	has no solution. If, on the contrary, there was  a solution $(\overline{x},\overline{t})\in\mathbb{R}^{n+1}$ such that
	$$
	\dfrac{1}{2}\overline{x}^\top A\overline{x}+\overline{t}a^\top \overline{x} + \overline{t}^2\gamma<0, \quad \max\left\{\dfrac{1}{2}\overline{x}^\top B\overline{x}+\overline{t}b^\top \overline{x}-\overline{t}^2\beta:(B,b)\in\mathcal{B}_b\right\}<0,
	$$
	$$
-\min\left\{\dfrac{1}{2}\overline{x}^\top B\overline{x}+\overline{t}b^\top \overline{x}-\overline{t}^2\alpha: (B,b)\in\mathcal{B}_b\right\}<0,
	$$
	we immediately reach a contradiction in case $\overline{t}\neq 0$. 
	So, suppose that $\overline{t}=0$. Then, the system \eqref{sis:hom1}-\eqref{sis:hom2} reduces to
	$$
	\dfrac{1}{2}\overline{x}^\top A\overline{x}<0, \quad	\max\left\{\dfrac{1}{2}\overline{x}^\top (B_1+\mu B_2)\overline{x}:\mu\in [\mu_1,\mu_2]\right\}<0,
	$$
	$$
	-\min\left\{\dfrac{1}{2}\overline{x}^\top (B_1+\mu B_2)\overline{x}: \mu\in [\mu_1,\mu_2]\right\}<0,
	$$
	which is impossible to hold. Thus, the claim is proved. \\
	On the other hand, observe that for every $(x,t)\in\mathbb{R}^n\times\mathbb{R}$, the minimum and maximum values in \eqref{sis:hom1}-\eqref{sis:hom2} are achieved in, at least, one of the extreme points of the rectangle $[\mu_1,\mu_2]\times[\delta_1,\delta_2]$, that is, in one of the elements  $(B_1+\mu_1B_2,b_1+\delta_1b_2)$, $(B_1+\mu_1B_2,b_1+\delta_2b_2)$, $(B_1+\mu_2B_2,b_1+\delta_1b_2)$ or  $(B_1+\mu_2B_2,b_1+\delta_2b_2)$. Then, the nonexistence of solution to the system \eqref{sis:hom1}-\eqref{sis:hom2} is equivalent to the nonexistence of $y\in\mathbb{R}^{n+1}$ solution to the system
	$$
	\dfrac{1}{2}y^\top H_0y<0
	$$
	$$
	\max\left\{\dfrac{1}{2} y^\top(W_{1_{\beta}}+\mu W_2)y
	:~\mu\in[\mu_1,\mu_2]\right\}<0,
	$$
	$$
	\max\left\{\dfrac{1}{2} y^\top(W_{2_{\beta}}+\mu W_2)y
	:~\mu\in[\mu_1,\mu_2]\right\}<0.
	$$
	$$
	-\min\left\{\dfrac{1}{2} y^\top(W_{1_{\alpha}}+\mu W_2)y:~\mu\in[\mu_1,\mu_2]\right\}
	<0,
	$$
	$$
	-\min\left\{\dfrac{1}{2} y^\top(W_{2_{\alpha}}+\mu W_2)y:~\mu\in[\mu_1,\mu_2]\right\}
	<0.
	$$
	
This means that	 $(0,0,0,0,0)\notin \Omega_W$.\\
\textbf{Step 2: A first use of a separation result.} 
Since $\Omega_W$ is convex, 
there exist  $(\lambda_0,\lambda_1,\lambda_2,\lambda_3,\lambda_4)\in\mathbb{R}^5_+\backslash\{0\}$ such that for all $y\in\mathbb{R}^{n+1}$
$$
\lambda_0 \left(\dfrac{1}{2}y^\top H_0y\right)+\lambda_1\left(\max_{\mu\in[\mu_1,\mu_2]}\dfrac{1}{2}y^\top (W_{1_{\beta}}+\mu W_2)y\right)+\lambda_2\left(\max_{\mu\in[\mu_1,\mu_2]}\dfrac{1}{2}y^\top (W_{2_{\beta}}+\mu W_2)y\right)+
$$
$$
\lambda_3\left(-\min_{\mu\in[\mu_1,\mu_2]}\dfrac{1}{2}y^\top (W_{1_{\alpha}}+\mu W_2)y\right)+\lambda_4\left(-\min_{\mu\in[\mu_1,\mu_2]}\dfrac{1}{2}y^\top (W_{2_{\alpha}}+\mu W_2)y\right)\geq 0.
$$
Since each of the minimum or maximum values are achieved in either $\mu_1$ or $\mu_2$, we obtain for all $y\in\mathbb{R}^{n+1}$,
$$
\lambda_0 \left(\dfrac{1}{2}y^\top H_0y\right)+\lambda_1\left(\max\left\{\dfrac{1}{2}y^\top (W_{1_{\beta}}+\mu_1 W_2)y, \dfrac{1}{2}y^\top (W_{1_{\beta}}+\mu_2 W_2)y\right\}\right)+
$$
$$
\lambda_2\left(\max\left\{\dfrac{1}{2}y^\top (W_{2_{\beta}}+\mu_1 W_2)y, \dfrac{1}{2}y^\top (W_{2_{\beta}}+\mu_2 W_2)y\right\}\right)+
$$
$$
\lambda_3\left(-\min\left\{\dfrac{1}{2}y^\top (W_{1_{\alpha}}+\mu_1 W_2)y, \dfrac{1}{2}y^\top (W_{1_{\alpha}}+\mu_2 W_2)y\right\}\right)+
$$
$$
\lambda_4\left(-\min\left\{\dfrac{1}{2}y^\top (W_{2_{\alpha}}+\mu_1 W_2)y, \dfrac{1}{2}y^\top (W_{2_{\alpha}}+\mu_2 W_2)y\right\}\right)\geq 0.
$$
Thus, there is no $y\in\mathbb{R}^{n+1}$ solution to the system:
$$
\dfrac{1}{2}y^\top \Big(\lambda_0 H_0+\lambda_1 (W_{1_{\beta}}+\mu_1 W_2)+\lambda_2 (W_{2_{\beta}}+\mu_1 W_2)-\lambda_3 (W_{1_{\beta}}+\mu_2 W_2)-\lambda_4 (W_{2_{\beta}}+\mu_2 W_2)\Big)y<0; 
$$
$$
\dfrac{1}{2}y^\top \Big(\lambda_0 H_0+\lambda_1 (W_{1_{\beta}}+\mu_2 W_2)+\lambda_2 (W_{2_{\beta}}+\mu_2 W_2)-\lambda_3 (W_{1_{\beta}}+\mu_1 W_2)-\lambda_4 (W_{2_{\beta}}+\mu_1 W_2)\Big)y<0. 
$$
This means that $(0,0)\not\in \Omega_2$, where 
$$
\Omega_2\doteq
\Big\{\Big(y^\top (\lambda_0 H_0+\lambda_1 
H_{1,1_{\beta}}+\lambda_2H_{1,2_{\beta}}
-\lambda_3 H_{2,1_{\beta}}-\lambda_4 H_{2,2_{\beta}})y,
$$
$$
y^\top(\lambda_0 H_0+\lambda_1 H_{2,1_{\beta}}
+\lambda_2 H_{2,2_{\beta}}
-\lambda_3 H_{1,1_{\beta}}
-\lambda_4 H_{1,2_{\beta}})
y\Big):~y\in\mathbb{R}^{n+1}\Big\},
$$
with 
	$$
W_{1_{\beta}}+\mu_1W_2=H_{1,1_{\beta}},\quad W_{1_{\beta}}+\mu_2W_2=H_{2,1_{\beta}},
$$
$$
W_{2_{\beta}}+\mu_1W_2=H_{1,2_{\beta}},\quad W_{2_{\beta}}+\mu_2W_2=H_{2,2_{\beta}},
$$
$$
W_{1_{\alpha}}+\mu_1W_2=H_{1,1_{\alpha}},\quad W_{1_{\alpha}}+\mu_2W_2=H_{2,1_{\alpha}},
$$
$$
W_{2_{\alpha}}+\mu_1W_2=H_{1,2_{\alpha}}\ \quad W_{2_{\alpha}}+\mu_2W_2=H_{2,2_{\alpha}}.
$$
\\
\textbf{Step 3: A second use of a separation result and conclusion.} 
By the Dines theorem, $\Omega_2$ is convex. Thus, there exists $(\xi_1,\xi_2)\in\mathbb{R}^2_+\backslash\{0\}$ such that for all $y\in\mathbb{R}^{n+1}$,
$$
\xi_1 y^\top \Big(\lambda_0 H_0+\lambda_1 H_{1,1_{\beta}}+\lambda_2H_{1,2_{\beta}}
-\lambda_3 H_{2,1_{\beta}}-\lambda_4 H_{2,2_{\beta}}\Big)y
$$
$$
+\xi_2y^\top\Big(\lambda_0 H_0+\lambda_1 H_{2,1_{\beta}}
+\lambda_2 H_{2,2_{\beta}}
-\lambda_3 H_{1,1_{\beta}}
-\lambda_4 H_{1,2_{\beta}}\Big)y\geq 0.
$$
In particular, for $y=(x,1)$ with $x\in\mathbb{R}^n$, and by setting 
$$\overline{\lambda}_0=\lambda_0(\xi_1+\xi_2),~~\overline{\lambda}_1=(\lambda_1+\lambda_2)(\xi_1+\xi_2),~~\overline{\lambda}_2=(\lambda_3+\lambda_4)(\xi_1+\xi_2),
$$ 
$$
\mu_{\beta}=\dfrac{\xi_1\mu_1+\xi_2\mu_2}{\xi_1+\xi_2},~
\mu_{\alpha}=\dfrac{\xi_1\mu_2+\xi_2\mu_1}{\xi_1+\xi_2},~\delta_\beta=\dfrac{\lambda_1\delta_1+\lambda_2\delta_2}{\lambda_1+\lambda_2},~{\rm and}~\delta_\alpha=\dfrac{\lambda_3\delta_2+\lambda_4\delta_1}{\lambda_3+\lambda_4},
$$
one gets  $(\overline{\lambda}_0,\overline{\lambda}_1,\overline{\lambda}_2)\in\mathbb{R}^3_+\backslash\{0\}$, $\mu_{\beta},\mu_{\alpha}\in[\mu_1,\mu_2]$, $\delta_{\beta},\delta_{\alpha}\in[\delta_1,\delta_2]$ and
$$
\overline{\lambda}_0\Big(\dfrac{1}{2}x^\top Ax+a^\top x+\gamma\Big) +\overline{\lambda}_1\Big(\dfrac{1}{2}x^\top (B_1+\mu_{\beta} B_2)x+(b_1+\delta_{\beta} b_2)^\top x-\beta\Big)+
$$
$$ 
\overline{\lambda}_2\left(\alpha-\Big(\dfrac{1}{2}x^\top (B_1+\mu_{\alpha} B_2)x+(b_1+\delta_{\alpha} b_2)^\top x \Big)\right)\geq 0 \quad \forall~ x\in\mathbb{R}^n,$$
where $\mu_{\alpha}+\mu_{\beta}=\mu_1+\mu_2$. This proves that $(b)$ holds.
\qed
\end{proof}


\begin{theorem} $($A robust S-lemma$)$\label{slema} 
	Let $A, B_1,B_2\in {\mathcal S}^n$, $a, b_1, b_2\in\mathbb{R}^n$ and $\gamma, \alpha, \beta$, $\mu_1, \mu_2,\delta_1,\delta_2\in\mathbb{R}$, with $\mu_1\leq\mu_2$, $\delta_1\leq\delta_2$ and $\alpha<\beta$. Assume that $\Omega_W$ is convex and that there exists 
	$x_0\in\mathbb{R}^n$ satisfying
 \begin{equation}\label{eq3}  \alpha<\dfrac{1}{2}x_0^\top (B_1+\mu B_2)x_0+(b_1+\delta b_2)^\top x_0 <\beta,~\forall~\mu\in [\mu_1,\mu_2], \forall~ \delta\in [\delta_1,\delta_2].
 \end{equation}
Then, the following two assertions are equivalent:
	\begin{itemize}
		\item[$(a)$] $\alpha\leq \dfrac{1}{2}x^\top (B_1+\mu B_2)x+(b_1+\delta b_2)^\top x \leq\beta$, $\forall~ \mu\in [\mu_1,\mu_2]$, $\forall ~\delta\in[\delta_1,\delta_2]$, $\Rightarrow \dfrac{1}{2}x^\top Ax+a^\top x+\gamma\geq 0$.
		\item[$(b)$] $\exists~(\lambda_1,\lambda_2)\in\mathbb{R}^2_+, \exists~\mu_{\alpha},\mu_{\beta}\in [\mu_1,\mu_2 ], \exists~\delta_{\alpha},\delta_{\beta}\in [\delta_1,\delta_2 ]:~ \forall~ x\in \mathbb{R}^n$
		$$
		\dfrac{1}{2}x^\top Ax+a^\top x+\gamma +\lambda_1\Big(	\dfrac{1}{2}x^\top (B_1+\mu_{\beta} B_2)x+(b_1+\delta_{\beta} b_2)^\top x-\beta\Big)+
		$$
		$$ \lambda_2\Big(\alpha-\Big(\dfrac{1}{2}x^\top (B_1+\mu_{\alpha} B_2)x+(b_1+\delta_{\alpha} b_2)^\top x \Big)\Big)\geq 0,
		$$
	\end{itemize}
	where $\mu_{\alpha}+\mu_{\beta}=\mu_1+\mu_2$.
\end{theorem}
\begin{proof} Clearly $(b)\Rightarrow(a)$.\\
Assume now that $(a)$ is satisfied. Then $(a)$ in Theorem \ref{teoralt} does not hold. Thus $(b)$ of the same theorem fulfills, but then $\lambda_0$ is strictly positive because of \eqref{eq3}, which implies the desired result.
\qed
\end{proof}

We are now ready to establish a characterization of optimality 
for the problem \eqref{prob:00}.

\begin{theorem}\label{minimonohomoalt}
	$($Characterizing robust optimality$)$   Let $A, B_1,B_2\in {\mathcal S}^n$, $a, b_1, b_2\in\mathbb{R}^n$ 
	and $\alpha, \beta, \mu_1, \mu_2,$ $\delta_1,\delta_2\in\mathbb{R}$, with $\mu_1\leq\mu_2$, $\delta_1\leq\delta_2$ and $\alpha<\beta$. 
	 Let 
$\overline{x}$ be feasible for problem \eqref{prob:00} and put $\gamma=-f(\overline{x})$. Assume that $\Omega_W$ is convex, and the 
Slater-type condition \eqref{eq3} holds. Then, $\overline{x}$ is optimal if, and only if there exist $(\lambda_1,\lambda_2)\in\mathbb{R}^2_+$, $\mu_{\alpha},\mu_{\beta}\in[\mu_1,\mu_2]$, $\delta_{\alpha}$, $\delta_{\beta}\in[\delta_1,\delta_2]$ 
	such that the following statements are satisfied:
	\begin{itemize}
		\item[$(a)$] $\Big(A+\lambda_1(B_1+\mu_{\beta}B_2)-\lambda_2(B_1+\mu_{\alpha}B_2)\Big)\overline{x}=-\Big(a+\lambda_1(b_1+\delta_\beta b_2)-\lambda_2(b_1+\delta_{\alpha}b_2)\Big);$
		\item[$(b)$] $\lambda_1\left(\dfrac{1}{2}\overline{x}^\top (B_1+\mu_{\beta}B_2)\overline{x}+(b_1+\delta_{\beta}b_2)^\top \overline{x}-\beta\right)=0;$\\
		$\lambda_2\left(\alpha-\Big(\dfrac{1}{2}\overline{x}^\top (B_1+\mu_{\alpha}B_2)\overline{x}+(b_1+\delta_{\alpha}b_2)^\top \overline{x}\Big)\right)=0;$
		\item[$(c)$] $A+\lambda_1(B_1+\mu_{\beta}B_2)-\lambda_2(B_1+\mu_{\alpha}B_2)\succeq 0$.
	\end{itemize}
\end{theorem}

\begin{proof} The necessary condition follows from Theorem \ref{slema} where $\gamma$ is substituted by $-f(\overline{x})$. Indeed, if $\overline{x}$ is optimal then $(a)$ in Theorem \ref{slema} holds, which means that $(b)$ of the same theorem is also satisfied. This finally implies the desired statements.\\
The sufficiency part is already standard since $(c)$ implies the convexity of the function
	$$
	h(x)\doteq\dfrac{1}{2}x^\top Ax+a^\top x-f(\overline{x}) +{\lambda}_1\Big(\dfrac{1}{2}x^\top (B_1+\mu_{\beta} B_2)x+(b_1+\delta_{\beta} b_2)^\top x-\beta\Big)+
	$$
	$$ 
	{\lambda}_2\left(\alpha-\Big(\dfrac{1}{2}x^\top (B_1+\mu_{\alpha} B_2)x+(b_1+\delta_{\alpha} b_2)^\top x \Big)\right),$$
	and $(a)$ and $(b)$ allow us to prove that $\overline{x}$ is in fact a solution to problem \eqref{prob:00}.
	\qed
\end{proof}

\section{Revisiting the case $\alpha=-\infty$}
\label{secc:uni}

We consider the problem:

\begin{equation}\label{uni0}
\begin{split}
\min &~ \dfrac{1}{2} x^\top Ax+a^\top x  \\
\text{s.t.}~ & \dfrac{1}{2}x^\top (B_1+\mu B_2)x+(b_1+\delta b_2)^\top\leq\beta,~ \forall~ \mu\in [\mu_1,\mu_2],\forall~\delta\in[\delta_1,\delta_2],
\end{split}
\end{equation}
where $A$, $B_1$, $B_2$ are real symmetric matrices, $\mu_1,\mu_2,\beta\in\mathbb{R}$ satisfying
$\mu_1< \mu_2$, $\delta_1<\delta_2$.
This problem was also discussed in \cite{jeya2013}, and actually this paper motivated our study.

By looking at carefully the proof of Theorem \ref{teoralt}, we immediately realize that in case there is no lower bound in the inequality constraint, all the terms where $\alpha$ appears, actually dissapear: they are superfluous. Hence, the set 
$\Omega_W$ reduces to
$$
\Omega^\beta_W\doteq\Big\{\Big(\dfrac{1}{2}y^\top H_0y
,\max_{\mu\in[\mu_1,\mu_2]}\dfrac{1}{2}y^\top (W_{1_{\beta}}+\mu W_2)y),\max_{\mu\in[\mu_1,\mu_2]}\dfrac{1}{2}y^\top (W_{2_{\beta}}+\mu W_2)y\Big):
$$
\begin{equation}
y\in\mathbb{R}^{n+1}\Big\}+\rm{int} \ \mathbb{R}^3_+.
\end{equation}
Thus, by Corollary 1 in \cite{fbp2022}, $\Omega_W^\beta$ is convex if the set
$$\Omega^\beta_\mu\doteq
\Big\{\Big(y^\top H_0y,y^\top(W_{1_{\beta}}+\mu_1 W_2)y,
y^\top(W_{2_{\beta}}+\mu_1 W_2)y, y^\top (W_{1_{\beta}}+\mu_2 W_2)y,
$$
$$
y^\top (W_{2_{\beta}}+\mu_2 W_2)y\Big):y\in\mathbb{R}^{n+1} \Big\}+\rm{int} \ \mathbb{R}^5_+.
$$

\begin{theorem}\label{minimocorregido} Let the data be as described above. Let $\overline{x}$ be feasible for  problem \eqref{uni0} and put $\gamma=-f(\overline{x})$. Assume that $\Omega_W^\beta$ is convex, and the 
Slater-type condition:  there exists $x_0\in\mathbb{R}^n$ such that 
\begin{equation}
\dfrac{1}{2}x_0^\top (B_1+\mu B_2)x_0+(b_1+\delta b_2)^\top x_0<\beta, \ \forall~ \mu\in[\mu_1,\mu_2], \ \forall ~\delta\in[\delta_1,\delta_2] 
\end{equation}
is satisfied. Then, $\overline{x}$ is optimal if, and only if there exist $\lambda\geq 0$, $\mu\in[\mu_1,\mu_2]$, $\delta\in[\delta_1,\delta_2]$ 
such that the following statements are satisfied:
\begin{itemize}
	\item[$(a)$] $\Big(A+\lambda(B_1+\mu B_2)\Big)\overline{x}=-\Big(a+\lambda(b_1+\delta b_2)\Big);$
	\item[$(b)$] $\lambda\left(\dfrac{1}{2}\overline{x}^\top (B_1+\mu B_2)\overline{x}+(b_1+\delta b_2)^\top \overline{x}-\beta\right)=0;$
	\item[$(c)$] $A+\lambda (B_1+\mu B_2)\succeq 0$.
\end{itemize}
\end{theorem}

We recall that
$$
H_0=\begin{pmatrix}A & a\\
a^\top  &  2\gamma \end{pmatrix},
$$
with $\gamma=-f(\overline{x})$ as in the previous theorem.
Set 
$$
H_1\doteq H_{1,1_{\beta}},~H_2\doteq H_{2,2_{\beta}}, H_3\doteq H_{1,2_{\beta}},~H_4\doteq H_{2,1_{\beta}}.
$$
The authors in \cite[Theorem 5.1]{jeya2013} proved the same result expressed in 
Theorem \ref{minimocorregido}  under the convexity of the set
\begin{equation}\label{setw00}
\Big\{\Big(y^\top H_0y,y^\top H_1y, y^\top H_2y\Big):~y\in\mathbb{R}^{n+1}\Big\}+
{\rm int}~\mathbb{R}_+^3;
\end{equation}
whereas ours requires the convexity of 
\begin{equation*}
\Omega_\mu^\beta=
\Big\{\Big(y^\top H_0y,y^\top H_1y, y^\top H_2y, y^\top H_3y, y^\top H_4y
\Big):~y\in\mathbb{R}^{n+1}\Big\}+{\rm int}~\mathbb{R}_+^5.
\end{equation*}
We believe that there is a gap in the proof of Theorem 5.1 in \cite{jeya2013}. More precisely, the authors assert in page 221 of the same paper that  the nonexistence of solution to the system (notice that our $\beta$ is $-\beta$ in \cite{jeya2013})
\begin{equation}\label{sis:hom10}
\dfrac{1}{2}x^\top Ax+ta^\top x + t^2\gamma<0,~ \max\left\{\dfrac{1}{2}x^\top Bx+tb^\top x-t^2\beta:(B,b)\in\mathcal{B}_b\right\}<0,
\end{equation}
implies that 
\begin{equation}\label{ecua00}
\dfrac{1}{2} y^\top  H_0 y<0\quad \text{and}\quad \forall~ \mu\in[\mu_1,\mu_2], \quad \dfrac{1}{2}y^\top (W_1+\mu W_2)y<0,
\end{equation}
has no solution. This is not necessarily true as Example \ref{ex100} below shows. We recall that
$$
W_1=\begin{pmatrix}B_1 & b_1+\dfrac{\delta_1\mu_2-\delta_2\mu_1}{\mu_2-\mu_1}b_2\\
(b_1+\dfrac{\delta_1\mu_2-\delta_2\mu_1}{\mu_2-\mu_1}b_2)^\top  &  -2\beta 
\end{pmatrix}
\quad \text{and}\quad 
$$
$$
W_{2}=\begin{pmatrix}B_2 & \dfrac{\delta_2-\delta_1}{\mu_2-\mu_1}b_2\\
\dfrac{\delta_2-\delta_1}{\mu_2-\mu_1}b_2^\top  &  0 \end{pmatrix}.
$$
\vskip0.3truecm
\begin{example} \label{ex100} Let  $n \geq 5$.  Taking $A=I_n$, $B_1=B_2=2I_n$,  $a=b_1=b_2 =s:=(1,1,\ldots,1) \in \mathbb{R}^n$,
$\mu_1=-1$,  $\mu_2=1$,  $\delta_1=-1, $  $\delta_2=1$ and $\beta =1$,  Problem \eqref{uni0} takes the form
\begin{equation}\label{example_j}
	\begin{split}
	\min ~~ & \dfrac{1}{2} \|x\|^2 + \sum_{i=1}^n x_i = \frac{1}{2} \|x+s\|^2 -\frac{n}{2}  \\
	\text{s.t.}~~ & g(x, \mu, \delta):=(1+\mu)\|x\|^2 + (1+\delta) \sum_{i=1}^n x_i -1 \leq 0~~\forall ~\mu,\delta\in[-1,1].
	\end{split}
	\end{equation}
	Let us define the functions
$$	
f(x):=\dfrac{1}{2} \|x\|^2 + \sum_{i=1}^n x_i = \frac{1}{2} \|x+s\|^2 -\frac{n}{2};
$$
$$
g(x, \mu, \delta):=(1+\mu)\|x\|^2 + (1+\delta) \sum_{i=1}^n x_i -1.
$$
Then, problem \eqref{example_j} is equivalent to 
	\begin{align}
	\min ~~& f(x)\notag\\ 
	\hbox{ s.t.}~~ &\|x\|^2 -\frac{1}{2} \leq 0; \label{eq1_example_j}\\
	& \sum_{i=1}^n x_i -\frac{1}{2} \leq 0;  \label{eq2_example_j}\\
	&\|x\|^2 +  \sum_{i=1}^n x_i -\frac{1}{2} \leq 0. \label{eq3_example_j}
	\end{align}
    Let $C$ be the set of constraints, that is, those $x$ satisfying 
    \eqref{eq1_example_j}-\eqref{eq3_example_j}. It is clear that $C$ is convex and compact. Thus, the unique solution to problem \eqref{example_j} is the projection
    of $-s$ on $C$, which is  $\overline{x} = -\dfrac{1}{\sqrt{2n}}s$. Hence,  
   $\overline{x}$ is a robust solution to problem \eqref{example_j}.  In  
	 this case, $\gamma = -f(\overline{x}) = \sqrt{\dfrac{n}{2}} -\dfrac{1}{4}$. By identifying the matrices involved in Theorem \ref{minimocorregido}, we get
	\[
	H_0=\begin{pmatrix}
	I_n & s \\
	s^\top & 2\gamma 
	\end{pmatrix};
	H_1= \begin{pmatrix}
	\overline{0} & 0 \\
	0^{\top} & -2 
	\end{pmatrix};
	H_2= \begin{pmatrix}
	4I_n & 2s \\
	2s^{\top} & -2
	\end{pmatrix};
	\]
	\[
	H_3=\begin{pmatrix} \overline{0} &2s\\
	2s^\top  &  2\beta \end{pmatrix};\quad H_4=\begin{pmatrix}4I_n  & 0\\
	0^\top  &  2\beta \end{pmatrix}.
	\]
	Since
	 \begin{equation} 
	-2H_0+ (-2-2\gamma) H_1+ H_2= \begin{pmatrix}
	2I_n & 0 \\
	0^{\top} & 2
	\end{pmatrix} \succ 0,
	\end{equation}
  applying \cite[Theorem 2.1 ]{Polyak} we have that 
\begin{equation}\label{setw:0}
\Big\{\Big(y^\top H_0y,y^\top H_1y, y^\top H_2y\Big):~y\in\mathbb{R}^{n+1}\Big\}
\end{equation}
is convex.  This means, according to \cite {jeya2013}, that  \eqref{example_j} is regular with respect to $\overline{x}$.
Also,  by taking $x_0=0$,   we have  $g(x_0, \mu, \delta)  <0$, for all  $\mu, \delta \in [-1, 1]$.  So,  we have all the conditions of   \cite[Theorem 5.1]{jeya2013} are satisfied. From the first part of the proof of \cite[Theorem 5.1]{jeya2013} we have that the following homogeneous system in $\mathbb{R}^{n+1}$, \eqref{sis:hom10} has no solution.\\
Coming back to our example, we obtain
\[W_1=\begin{pmatrix}
	2I_n & s \\
	s^{\top} & -2 
	\end{pmatrix} \hbox{ and }
	W_2=\begin{pmatrix}
	2I_n & s \\
	s^{\top} & 0 
	\end{pmatrix}.
	\]
Then, the following homogeneous system in $\mathbb{R}^{n+1}$ (see \eqref{ecua00})
\begin{equation*} 
\frac{1}{2}\begin{pmatrix} x\\t\end{pmatrix}^{\top} H_0\begin{pmatrix} x\\t\end{pmatrix}  <0 \hbox{ and }  
\forall ~\mu \in [-1, 1],~\frac{1}{2}\begin{pmatrix} x\\t\end{pmatrix}^{\top}(W_1+\mu W_2)\begin{pmatrix} x\\t\end{pmatrix}<0,
\end{equation*}
becomes
\begin{equation}\label{ecua01}
\frac{1}{2}\begin{pmatrix} x\\t\end{pmatrix}^{\top}\begin{pmatrix}
	I_n & s \\
	s^{\top} & 2\gamma 
	\end{pmatrix}\begin{pmatrix} x\\t\end{pmatrix}  <0 \hbox{ and }  
\end{equation}
\begin{equation}\label{ecua02}
\forall~ \mu \in [-1, 1],~\frac{1}{2}\begin{pmatrix} x\\t\end{pmatrix}^{\top}\begin{pmatrix}
	2(1+\mu)I_n & \; & (1+\mu)s \\
	(1+\mu)s^{\top} & \;& -2
	\end{pmatrix}\begin{pmatrix} x\\t\end{pmatrix}<0.
\end{equation}
We will see that such a system admits a solution, contradicting the assertion made in page 221 of \cite{jeya2013} about the nonexistence of solution to the same system. Indeed, for $(-s,1) \in \mathbb{R}^{n+1}$, \eqref{ecua01} reduces to 
\[
\frac{1}{2}\begin{pmatrix} -s\\1\end{pmatrix}^{\top}\begin{pmatrix}
	I_n & s \\
	s^{\top} & 2\gamma 
	\end{pmatrix} \begin{pmatrix} -s\\1\end{pmatrix} = \frac{1}{2}(-n+2\gamma )= \frac{1}{2}\big(-n+  \sqrt{\frac{n}{2}} -\frac{1}{4}\big) < 0;
\]
whereas \eqref{ecua02} becomes: for all $\mu \in [-1, 1]$,
	 \[
	\frac{1}{2}\begin{pmatrix} -s\\1\end{pmatrix}^{\top}\begin{pmatrix}
	2(1+\mu)I_n & \;& (1+\mu)s \\\
	(1+\mu)s^{\top} & \; & -2
	\end{pmatrix}\begin{pmatrix} -s\\1\end{pmatrix}= 	\frac{1}{2}(n(1+\mu)-(1+\mu)n -2) =-1 <0.
	\]
This proves our claim. 

Observe that taking  in Corollary \ref{cor_R_5}: $m=4$, $A_0=I_n$, $\rho_1=\rho_3=4$,  $\rho_2=\rho_4=0$,
$a_0= s, a_1=a_3=0$,  $a_2=a_4=2s$,  $c_0= \gamma$, $c_1=c_2= -1$, $c_3=c_4= \beta$ we obtain that
$$
\Big\{\Big(y^\top H_0y,y^\top H_1y, y^\top H_2y, y^\top H_3y, y^\top H_4y\Big):~y\in\mathbb{R}^{n+1}\Big\} \subset \mathbb{R}^5
$$
is a convex set, which is required in our Theorem \ref{minimocorregido}, providing the characterization of robust optimality for our example. 
\qed\end{example}

Clearly, the previous example only shows there is a gap in the proof of Theorem 5.1 in \cite{jeya2013}. We were unable to construct a real counterexample to that result under the convexity either of the set given in \eqref{setw00} or in \eqref{setw:0}.

\end{document}